\documentclass[12pt]{amsart}
\usepackage{amsmath,amsthm}

\def\pmatrix{\left(\begin{matrix}}
\def\endpmatrix{\end{matrix}\right)}

\def\P{{\mathbb P}}
\def\Z{{\mathbb Z}}
\def\C{{\mathbb C}}

\def\cal{\mathcal}

\def\dd{{\rm d}}

\def\diag{\operatorname{diag}}
\def\de{\delta}

\def\p{\partial}

\def\t{\theta}
\def\s{\sigma}
\def\T{\Theta}
\def\e{\varepsilon}

\def\a{\alpha}
\def\b{\beta}

\def\A{{\mathcal A}}

\def\J{{\mathcal J}}
\def\H{{\mathcal H}}

\def\tt#1#2{{\t\left[\begin{matrix}{#1}\\ {#2}\end{matrix}\right]}}

\theoremstyle{plain}
\newtheorem{thm}{Theorem}
\newtheorem{lm}[thm]{Lemma}
\newtheorem{prop}[thm]{Proposition}
\newtheorem{cor}[thm]{Corollary}

\theoremstyle{definition}

\newtheorem{rem}[thm]{Remark}

\font\xxx=cmti8 \font\yyy=msbm8
\title{Gradients of odd theta functions}
\author{Samuel Grushevsky\and Riccardo Salvati Manni}
\email{sam@math.princeton.edu\and salvati@mat.uniroma1.it}
\address{Mathematics Department, Princeton University, Fine Hall,
Washington Road, Princeton, NJ 08544, USA\and Dipartimento di
Matematica, Universit\`a di Roma, Piazzale Aldo Moro, 2, I-00185
Roma, Italy}
\date{May 11, 2003.}
\thanks{First author partially supported by NSF Mathematical Sciences
Postdoctoral Research Fellowship} 

\begin{document}

\begin{abstract}
We show that a generic principally polarized abelian variety is
uniquely determined by its gradient theta-hyperplanes, the
non-projectivized version of those studied in \cite{cap1},
\cite{cap2}, \cite{cap4}, which in a sense are a generalization to
ppavs of bitangents of plane curves. More precisely, we show that,
generically, the set of gradients of all odd theta functions at
the point zero uniquely determines a ppav with level (4,8)
structure. We also show that our map is an immersion of the moduli
space of ppavs.
\end{abstract}
\maketitle

\section{Definitions and notations}
We denote by $\H_g$ the {\it Siegel upper half-space} --- the
space of complex symmetric $g\times g$ matrices with positive
definite imaginary part. An element $\tau\in\H_g$ is called a {\it
period matrix}, and defines the complex abelian variety
$X_\tau:=\C^g/\Z^g+\tau \Z^g$. The group $\Gamma_g:={\rm
Sp}(2g,\Z)$ acts on $\H_g$ by automorphisms: for $\gamma:=\pmatrix
a&b\\ c&d\endpmatrix\in{\rm Sp}(2g,\Z)$ the action is
$\gamma\tau:=(a\tau+b)(c\tau+d)^{-1}$. The quotient of $\H_g$ by
the action of the symplectic group is the moduli space of
principally polarized abelian varieties (ppavs): $\A_g:=\H_g/{\rm
Sp}(2g,\Z)$. A ppav is called irreducible if it is not a direct
product of two lower-dimensional ppavs, i.e. if its period matrix
$\tau$ is not conjugate by the action of $\Gamma_g$ to a matrix
that splits as $\tau_1\oplus\tau_2$ for two lower-dimensional
period matrices. For us the case $g=1$ is special and in the
following we will always assume $g>1$.

We define the {\it level} subgroups of the symplectic group to be
$$
\Gamma_g(n):=\left\lbrace\gamma=\pmatrix a&b\\ c&d\endpmatrix
\in\Gamma_g\, |\, \gamma\equiv\pmatrix 1&0\\
0&1\endpmatrix\ {\rm mod}\ n\right\rbrace
$$
$$
\Gamma_g(n,2n):=\left\lbrace\gamma\in\Gamma_g(n)\, |\, {\rm diag}
(a^tb)\equiv{\rm diag}(c^td)\equiv0\ {\rm mod}\ 2n\right\rbrace.
$$
The corresponding {\it level moduli spaces of ppavs} are denoted
$\A_g^n$ and $\A_g^{n,2n}$, respectively.

A function $F:\H_g\to\C$ is called a {\it modular form of weight $k$
with respect to $\Gamma\subset\Gamma_g$} if
$$
F(\gamma\tau)=\det(c\tau+d)^kF(\tau),\quad \forall \gamma=
\pmatrix a&b\\ c&d\endpmatrix\in\Gamma,\ \forall \tau\in\H_g
$$

More generally, let $\rho:{\rm GL}(g,\C)\to\operatorname{End} V$ be some
representation. Then a map $F:\H_g\to V$ is called a {\it $\rho$- or
$V$-valued modular form}, or simply a {\it vector-valued modular form},
if the choice of $\rho$ is clear, with respect to
$\Gamma\subset\Gamma_g$ if
$$
F(\gamma\tau)=\rho(c\tau+d)F(\tau),\quad \forall \gamma=
\pmatrix a&b\\ c&d\endpmatrix\in\Gamma,\
\forall \tau\in\H_g.
$$

For $\e,\de\in \Z_2^g$, thought of as vectors of zeros and ones,
and $z\in \C^g$ we define the {\it theta function with
characteristic $[\e,\de]$} to be
$$
\tt\e\de(\tau,z):=\sum\limits_{m\in\Z^g} \exp \pi i \left[\left(
m+\frac{\e}{2},\tau(m+\frac{\e}{2})\right)+2\left(m+\frac{\e}{2},z+
\frac{\de}{2}\right)\right].
$$
A {\it characteristic} $[\e,\de]$ is called {\it even} or {\it
odd} depending on whether the scalar product $\e\cdot\de\in\Z_2$
is zero or one, and the corresponding   theta function is even or
odd in $z$, respectively. The number of even (resp. odd)
characteristics is $2^{g-1}(2^g+1)$ (resp. $2^{g-1}(2^g-1)$).

For $\e\in\Z_2^g$ we also define the {\it second order theta
function} with characteristic $\e$ to be
$$
\T[\e](\tau,z):=\tt{\e}{0}(2\tau,2z).
$$
The group $\Gamma_g$ acts on the set of characteristics  as
follows:
$$
\gamma \pmatrix\e\\ \de\endpmatrix  :=\pmatrix d &  -c \\
-b  & a \endpmatrix \pmatrix \e\\ \de \endpmatrix +\frac{1}{2}
\pmatrix \diag(cd^t)\cr \diag(ab^t) \endpmatrix,
$$
where the resulting characteristics is taken modulo 2. This action
is not transitive, in fact the parity of the characteristics is an
invariant. The transformation law for theta functions under the
action of the symplectic group is (see \cite{ig1}):
\begin{equation}
\label{translaw} \theta\left[\gamma \pmatrix\e\\ \de \endpmatrix
\right](\gamma\tau,(c\tau+d)^{-t}
z)=\phi(\e,\de,\gamma,\tau,z)\det (c\tau+d)^{1/2}\tt\e\de(\tau,z),
\end{equation}
where $\phi$ is some complicated explicit function. It is further
known (see \cite{ig1}, \cite{sm1}) that $\phi|_{z=0} $ does not
depend on $\tau$, and that for $\gamma\in\Gamma_g(4,8)$ we have
$\phi|_{z=0}=1$, while $\gamma\in\Gamma_g(4,8)$ acts trivially on
the characteristics $[\e,\de]$. Thus the values of theta functions
at $z=0$, called {\it theta constants}, are modular forms of
weight one half with respect to $\Gamma_g(4,8)$. Similarly it is
known that the theta constants of second order are modular forms
of weight one half with respect to $\Gamma_g(2,4)$. The action of
$\Gamma_g(2)/ \Gamma_g(4,8)$ on the set of theta constants with
characteristics is by certain characters whose values are fourth
roots of the unity, and is well understood --- see \cite{sm1}. The
action of $\Gamma_g/\Gamma_g(2)$ on the set of characteristics is
by permutations.

All odd theta constants with characteristics vanish identically,
as the corresponding theta functions are odd functions of $z$, and
thus there are $2^{g-1}(2^g+1)$ non-trivial theta constants with
characteristics, and $2^g$ theta constants of the second order.

Differentiating the theta transformation law above with respect to
$z_i$ and then evaluating at $z=0$, we see that for $\gamma\in
\Gamma_g(4,8)$ and $[\e,\de]$ odd
$$
\frac{\p}{\p z_i}\tt\e\de(\tau,z)|_{z=0}=\det(c\tau+d)^{1/2}
\sum\limits_j (c\tau+d)_{ij}\frac{\p}{\p
z_j}\tt\e\de(\gamma\tau,z)|_{z=0};
$$
in other words the gradient vector $\lbrace\frac{\p}{\p
z_i}\tt\e\de(\tau,0)\rbrace_{{\rm all}\ i}$ is a $\C^g$-valued
modular form with respect to $\Gamma_g(4,8)$ under the
representation $\rho(M):=(\det M)^{1/2}\cdot M$.

The set of all even theta constants with characteristics defines
the map
$$
\P Th:\A_g^{4,8}\to \P^{2^{ g-1}(2^g+1)-1},\quad \P Th(\tau):=\lbrace
\tt\e\de(\tau,0) \rbrace_{{\rm all\ even}\ [\e,\de]}
$$
Theta constants of the second order similarly define the map
$$
\P Th_2:\A_g^{2,4}\to \P^{2^g-1},\quad \P Th_2(\tau):=\lbrace\T
[\e](\tau,0) \rbrace_{{\rm all}\ \e}.
$$
Considering the set of gradients of all odd theta functions at
zero gives the map
$$
grTh:\H_g\to (\C^g)^{\times 2^{g-1}(2^g-1)}\qquad grTh(\tau):=\left\lbrace
\vec{\rm grad}_z\tt\e\de\right\rbrace_{{\rm all\ odd}\ [\e,\de]},
$$
which due to modular properties descends to the quotient map
$$
\P grTh:\A_g^{4,8}\to  (\C^g)^{\times 2^{g-1}(2^g-1)}/\rho({\rm GL}(g,\C)),
$$
where ${\rm GL}(g,\C)$ acts simultaneously on all $\C^g$'s in the
product by $\rho$.

Because of Lefschetz theorem for abelian varieties for any $\tau$
the rank of the $2^{g-1}(2^g-1)\times g$ matrix of derivatives
$\frac{\p}{\p_{z_i}}\tt \e\de(\tau,0)$ is always $g$ (see
\cite{sm3}). Thus if we think of this matrix as a $g$-tuple of
vectors in $\C^{2^{g-1}(2^g-1)}$, it is always non-degenerate.
Thus the image of $\P grTh$ in fact lies in the grassmannian,
\begin{equation}
\P grTh:\A_g^{4,8}\to{\rm Gr}_\C(g,\,2^{g-1}(2^g-1))
\end{equation}
of $g$-dimensional subspaces in $\C^{2^{g-1}(2^g-1)}$. The
Pl\"ucker's coordinates of this map are modular forms of weight
$\frac{g}{2}+1$ and have been extensively studied --- see
\cite{fr885}, \cite{fay}, \cite{ig2}, \cite{sm3}.

It is known that the map $\P Th$ is an  embedding  --- see
\cite{ig1} and references therein. In \cite{sm1} it is shown that
$\P Th_2$ is also injective. However, there appears to be a small
gap in the proof there, so at the moment we can only say that the
map is injective for $g \leq 3$ and generically injective for
$g\ge 4$. In this work we will avoid using this result except for
the case $g=2$. We remark that, because of the transformation
formula of theta functions, all the  above maps are
$\Gamma_g$-equivariant. Here we prove the following properties of
$\P grTh$:

\begin{thm}
For $g\geq 3$ the map $\P grTh$ is generically injective on
$\A_g^{4,8}$. For genus 2 the map is finite of degree 16.
\end{thm}

In the course of the proof we give explicitly an open set in
$\A_g^{4,8}$ where the map is injective, and obtain, denoting by
$\J_g\subset\A_g$ the locus of Jacobians,
\begin{cor}
For $g\geq 3$ the map $\P grTh$ is also generically injective on
$\J_g^{4,8}$.
\end{cor}

We will also show
\begin{thm}
For $g\geq 2$ the map $\P grTh$ is injective on tangent spaces.
\end{thm}

\begin{rem}
The map $\P grTh$ is closely related to the indicated in
\cite{cap2} generalization to abelian varieties of the map
obtained by sending a Jacobian to the collection of the
hyperplanes tangent to the canonical curve at $g-1$ points,
considered in \cite{cap1}, \cite{cap2}, \cite{cap4}. This map
itself is in fact the generalization of the map of a plane quartic
to the set of its bitangents. It was shown by Aronhold in
\cite{aronh} that some subset of bitangents (i.e. their
directions), now called an Aronhold system, together with the
points of their tangency serves to recover the curve, but until
\cite{cap2} it was not known that the bitangent directions
themselves determine a (generic) curve.

Indeed, suppose we have a hyperplane tangent to the canonical
curve at $g-1$ points. The reduced divisor that it cuts on the
curve is then the square root of the canonical, and thus a theta
characteristic. However, it is an effective theta characteristic,
which for a generic curve means it is an odd theta characteristic.
On the other hand, for each odd theta characteristics we get a
$(g-1)$-tangent hyperplane to the canonical curve, and
analytically it is clear what the direction of this hyperplane is:
it is given by the gradient of the corresponding odd theta
function with characteristics at zero.

In our terms the map to $(g-1)$-tangent hyperplanes is thus the
map
$$
gr\P Th:\J_g\to S^{2^{g-1}(2^g-1)}(\P^{g-1})/{\rm PGL}(g,\C)
$$
of the locus of Jacobians $\J_g$, obtained by sending $\tau$ to
the set of gradients of all odd theta constants, but {\it each}
projectivized independently, and considered as a point in
$\P^{g-1}$ and not in $\C^g$ (we have the symmetric power of
$\P^{g-1}$ here instead of the direct power because we are
forgetting the level structure, and thus the characteristics may
be permuted).

This map $gr\P Th$ considered by Caporaso and Sernesi is quite
different from our $\P grTh$. As explained in \cite{cap2} and
\cite{cap4}, $gr\P Th$ is not defined on all of $\J_g$ or $\A_g$
--- the problem occurs for those $\tau$ for which one of the
components of $grTh(\tau)$ is equal to zero. Where it is defined,
it factors through $\P grTh$:
$$
\begin{matrix}
&\A_g&\mathop{\longrightarrow}\limits^{\P grTh}\quad&{\rm
Gr}_\C(g,\,2^{g-1}(2^g-1))\ \\
&&\quad\searrow \hbox{\xxx gr\yyy P\xxx Th}\quad&\downarrow\pi\\
&&&(\P^{g-1})^{2^{g-1}(2^g-1)}/{\rm PGL}(g,\C)\ ,
\end{matrix}
$$
where $\pi$ is the natural projection, which is not always
defined.

In \cite{cap2} and \cite{cap4} it is shown that when restricted to
its domain on the Jacobian locus, $gr\P Th$ is generically
injective. While our corollary 2 states that $\P grTh|_{\J_g}$ is
generically injective, it does not serve to reproduce the result
of Caporaso and Sernesi, as we do not handle the projectivization
map $\pi$, which might collapse images different points.
\end{rem}

In the following, we will often omit the arguments $\tau$ and
$z=0$ for theta functions and their derivatives, and will  write
$\p_{z_i}$ instead of $\frac{\p}{\p z_i}$. The letters of the
Latin alphabet will denote coordinates for vectors in $\C^g$, i.e.
will range from 1 to $g$. The letters of the Greek alphabet denote
characteristics for $g$-dimensional theta functions, i.e. lie in
$\Z_2^g$.

\section{$\t$'s and $\T$'s}
The fundamental relation between theta functions with
characteristics and theta functions of the second order comes from
the fact that the squares of theta functions with characteristics
are sections of the bundle on the abelian variety for which the
theta functions of thesecond order form a basis for the space of
sections. The relationship and others more general are special
cases of Riemann's addition theorem for theta functions (see, for
example, \cite{ig1}):
\begin{equation}
\label{tT}
\begin{matrix}
\tt\a\b(2\tau,2z)\tt{\a+\e}\b(2\tau,2x)\\
=\frac{1}{2^g}\sum\limits_{\s\in\Z_2^g}(-
1)^{\a\cdot\s}\tt\e{\b+\s}(\tau,z+x)\tt\e\s (\tau,z-x),
\end{matrix}
\end{equation}
which is valid for all $\tau$, $z$, $x$ and $\a,\b,\e$. Let us
denote
\begin{equation}
\label{defC}
\begin{matrix}
C_{ij\,\e\de}^{\b}(\tau):=\p_{z_i}\tt\e{\b+\de}(\tau,0)\p_{z_j}
\tt\e\de(\tau,0)\\
+\p_{z_j}\tt\e{\b+\de}(\tau,0)\p_{z_i}\tt\e\de(\tau,0),
\end{matrix}
\end{equation}
which is zero unless $[\e,\b+\de]$ and $[\e,\de]$ are odd, and
\begin{equation}
\begin{matrix}
A_{ij\,\e\de}^{\b}(\tau):=\p_{z_i}\p_{z_j}\tt\de\b(2\tau,0)\tt\e\b(2\tau,0)-
\\
\tt\de\b(2\tau,0)\p_{z_i}\p_{z_j}\tt\e\b(2\tau,0),
\end{matrix}
\end{equation}
which is zero unless $[\de,\b]$ and $[\e,\b]$ are both even
characteristics. Note also that by the heat equation we have
$\p_{z_i}\p_{z_j}\t=\p_{ \tau_{ij}}\t$ up to a constant that is
not important to us.

We then have the following relation
\begin{lm}
If $[\e,\de]$ and $[\e,\b+\de]$ are odd characteristics,
\begin{equation}
C_{ij\,\e\de}^{\b}= \frac{1}{2}\sum\limits_{\a\in\Z_2^g}(-1)^{
\a\cdot\de }A_{ij\, {\e+\a} \a}^{\b}.
\end{equation}
\end{lm}
\begin{proof}
Let us take the sum of the equations (\ref{tT}) for different
$\a$, each with coefficient $(-1)^{\a\cdot\de}$, where $\de$ is
some characteristic. We get
\begin{equation}
\label{tT2}
\begin{matrix}
\sum\limits_{\a\in\Z_2^g}
(-1)^{\a\cdot\de}\tt\a\b(2\tau,2z)\tt{\a+\e}\b(2\tau,0)\\
={1\over 2^g}\sum\limits_{\a,\s\in\Z_2^g}(-1)^{
\a\cdot(\s+\de)}\tt\e{\b+\s}(\tau,z)\tt\e\s (\tau,z)\\
=\tt\e{\b+\de}(\tau,z)\tt\e\de (\tau,z).
\end{matrix}
\end{equation}
Now differentiate this relation twice with respect to $z_i$ and
$z_j$ and then evaluate at $z=0$ to prove the lemma.
\end{proof}
Moreover, the expression of $C$'s in terms of $A$'s is invertible:

\begin{lm}
\begin{equation}
\label{Adef}
A_{ij\,{\a+\e}\a}^{\b}=\frac{1}{2^{g-1}}\sum_{\lbrace\s|[\e,\s]
{\rm\ odd}\rbrace}(-1)^{\a\cdot\s}C_{ij\,\e\s}^\b
\end{equation}
\end{lm}
\begin{proof}
In (\ref{tT}), we assume that $[\a,\b] $ and $[\a+\e,\b]$ are even
characteristics. Differentiating, we get
$$
\p_{z_i}\p_{z_j}\tt\a\b(2\tau,0) \tt{\a+\e}\b(2\tau,0)
$$
$$
=\frac{1}{2^{g}}\p_{z_i}\p_{z_j}\left.\left(\sum\limits_{\s\in\Z_2^g}(-1)^{\a\cdot\s}
\tt\e{\b+\s}( \tau, z) \tt\e\s (\tau,z)\right)\right|_{z=0}.
$$
Similarly, switching $\a$ and $\a+\e$ we have
$$
\tt\a \b ( 2\tau,
0)\p_{z_i}\p_{z_j} \tt{\a+\e}\b (2\tau,0)
$$
$$
=\frac{1}{2^{g}}\p_{z_i}\p_{z_j}\left.\left(\sum\limits_{\s\in\Z_2^g}(-
1)^{(\a+\e)\cdot\s} \tt\e{\b+\s}( \tau, z) \tt\e\s
(\tau,z))\right)\right|_{z=0}.
$$
Subtracting and computing separately for the cases of $[\e,\s]$
odd and even, we get the statement of the lemma.
\end{proof}

\section{Recovering $\P Th_2(\tau)$ from $\P grTh(\tau)$}
\begin{prop}
The following identity holds for all $i,j,\e,\de,\s,\b$, and all
$\tau\in\H_g$:
\label{propline}
\begin{equation}
\label{syseq}
A_{ij,\e\de}^{\b}(\tau)\tt\s\b(2\tau)+A_{ij,\de\s}^{\b}(\tau)\tt\e\b
(2\tau)+A_{ij,\s\e}^{\b}(\tau)\tt\de\b(2\tau)=0.
\end{equation}
\end{prop}
\begin{proof}
From the definition of $A$'s it follows that the above expression
is the determinant of the matrix
$\left(\begin{matrix}
\tt\e\b&\tt\e\b&\p_{\tau_{ij}}\tt\e\b\\
\tt\de\b&\tt\de\b&\p_{\tau_{ij}}\tt\de\b\\
\tt\s\b&\tt\s\b&\p_{\tau_{ij}}\tt\s\b
\end{matrix}\right)$.
\end{proof}

\begin{prop}
\label{propminor}
For any $i,j,I,J,\e,\de,\s,\tau$ we have
\begin{equation}
\label{betaalpha}
\begin{matrix}
\left(A_{ij,\e\de}^\b(\tau) A_{IJ,\s\de}^\b(\tau)-A_{IJ,\e\de}^\b(\tau)
A_{ij,\s\de}^\b(\tau)\right)\tt\e\b(2\tau)\\
=\left(A_{ij,\de\e}^\b(\tau) A_{IJ,\s\e}^\b(\tau)-A_{IJ,\de\e}^\b(\tau)
A_{ij,\s\e}^\b(\tau)\right)\tt\de\b(2\tau).
\end{matrix}
\end{equation}
\end{prop}
\begin{proof}
The identity is straightforward if we substitute the definition of
$A$'s in terms of $\t$'s and their derivatives. Indeed
\begin{equation}
\label{minors}
A_{ij,\e\de}^\b A_{IJ,\s\de}^\b-A_{IJ,\e\de}^\b A_{ij,\s\de}^\b=
\tt\de\b\det\pmatrix
\tt\e\b & \p_{\tau_{ij}}\tt\e\b & \p_{\tau_{IJ}}\tt\e\b\\
\tt\de\b & \p_{\tau_{ij}}\tt\de\b & \p_{\tau_{IJ}}\tt\de\b\\
\tt\s\b & \p_{\tau_{ij}}\tt\s\b & \p_{\tau_{IJ}}\tt\s\b
\endpmatrix
\end{equation}
evaluated at $2\tau$. Another proof would be to use the previous
proposition for $\e,\de,\s$ with $i,j$ and then with $I,J$ to
express $\tt\s\b$ in two different ways, and then equate these two
expressions.
\end{proof}

Let us consider the case $\b=0$ for the propositions above. Then
the $A_{ij,\e\de}^{0}$ are $2\times 2$ minors of the $2^g\times
\left(\frac{g(g+1)}{2}+1\right)$ matrix $M$ with columns
$(\T[\e],\lbrace\p_{\tau_{ij}}\T[\e]\rbrace_{ i\le j})_{\rm all\
i,j}$ (i.e. in each column there are the values either of $\T[\e]$
or of its derivative, for all $\e\in\Z_2^g$). This has maximal
rank, equal to $\frac{g(g+1)}{2}+1$, for irreducible abelian
varieties
--- see \cite{sasaki} and \cite{sm1}. Inductively for reducible
abelian varieties we see that for $g\ge 2$ the rank of $M$ is at
least three.

\begin{prop}
For $\b=0$ the system of equations (\ref{propline}) has a unique
projective solution for $\lbrace\T[\e](\tau)\rbrace_{{\rm all}\
\e}$ for fixed $A$'s. Since $A$'s are expressible in terms of
$C$'s, which are combinations of gradients of odd theta functions,
this means that $\P Th_2(\tau)$ is determined uniquely by $\P
grTh(\tau)$.
\label{getTh2}
\end{prop}
\begin{proof}
This is a consequence of the fact that the matrix $M$ above has at
least three, and that the solutions of the system (\ref{propline})
are invariant under the action of ${\rm GL}(g,\C)$. Basically we
need to show that the system has maximal rank. Suppose we are
given $\P grTh(\tau)$, i.e. all $A^0_{ij,\e\de}$'s. Let us pick a
representative $\tau\in\H_g$ and think of $grTh(\tau)$ --- we will
deal with the action of ${\rm GL}(g,\C)$ later.

Since the matrix $M$ described above has rank at least 3, we can
pick a non-degenerate $3\times 3$ minor in it. For irreducible
$\tau$ the matrix $M$ is of maximal rank, and thus this minor can
be chosen to contain the first column. From the fact that the
theta constants of reducible abelian varieties are products of
lower-dimensional theta constants it follows that such a choice is
also possible for reducible abelian varieties.

Suppose now that this non-degenerate minor is
$$
\det\pmatrix
\T[\e] &\p_{\tau_{ij}}\T[\e] &\p_{\tau_{IJ}}\T[\e]\\
\T[\de] &\p_{\tau_{ij}}\T[\de] &\p_{\tau_{IJ}}\T[\de]\\
\T[\s] &\p_{\tau_{ij}}\T[\s] &\p_{\tau_{IJ}}\T[\s]
\endpmatrix\ne 0.
$$
Then at least one of $\T[\e],\T[\de],\T[\s]$ must also be non-zero
--- by renaming let it be $\T[\e]$. Then the combination of $A$'s
in the right-hand-side of the formula (\ref{minors}) is non-zero,
and thus we can use proposition \ref{propminor} to express
$\frac{\T[\de]}{\T[\e]}$ and $\frac{\T[\s]}{\T[\e]}$ in terms of
$A$'s, i.e. in terms of $grTh(\tau)$.

Furthermore, since $\T[\e]$ and the $3\times 3$ minor above are
non-zero, some $2\times 2$ subminor containing $\T[\e]$ must also
be non-zero. By renaming let it be
$$
\det\pmatrix
\T[\e] &\p_{\tau_{ij}}\T[\e]\\
\T[\de] &\p_{\tau_{ij}}\T[\de]\\
\endpmatrix=A_{ij,\e\de}^0\ne 0.
$$
Then we can use proposition \ref{propline} with $\e,\de,i,j$ and
any $\s$ to express all $\frac{\T[\s]}{\T[\e]}$ in terms of $A$'s.
Thus from $grTh(\tau)$ we can recover $\P Th_2(\tau)$ uniquely.

Now we have to deal with the action of ${\rm GL}(g,\C)$ to finish
the proof. However, the system (\ref{syseq}) is acted upon by the
adjoint action of ${\rm GL}(g,\C)$ (if we consider each
$A_{ij\,\e\de}^0$ as a matrix labeled by $i,j$), and thus
transformed into an equivalent system. This equivalent system will
have the same solutions, and as we have shown the solution to be
unique, it has the same solution.
\end{proof}

\section{Generic injectivity of $\P grTh$}
If it were known that $\P Th_2$ is injective, we would be already
done, and could conclude the injectivity of $\P grTh$ at level
(2,4). But, as this is not yet known, we need to do extra work.

\begin{rem}
The genus $g=2$ case is rather special: in this case it is known
that $\P Th_2$ is injective; however, there are only six odd
characteristics, while
$$
\vert \Gamma_2(2,4)/\Gamma_2(4,8)\vert=2^{10},
$$
Thus in genus two the map $\P grTh(\tau)$ factors over some
subgroup $\Gamma$ such that $\Gamma_2(4,8)\subset \Gamma\subset
\Gamma_2(2,4)$, and the injectivity of $\P grTh$ holds only on
$\A_2^\Gamma$.
\end{rem}

In general, to recover $\P Th(\tau)$ from $\P Th_2(\tau)$, we just
need to know which sign to choose for each $\tt\e\de$, since their
squares are already expressible in terms of theta constants of the
second order by (\ref{tT2}).

\smallskip
To get more control over the signs, let us consider the case
$\b\neq 0$ of the equations (\ref{tT}). Then $A_{ij,\e\de}^{\b}$
are $2\times 2$ minors of the $2^g\times
\left(\frac{g(g+1)}{2}+1\right)$ matrix $M^{\b}$ consisting of
$(\tt\a\b(2\tau,0),\lbrace\p_{\tau_{ij}} \tt \a\b (2\tau,0)
\rbrace_{i\le j})_{{\rm all\ }\a}$, which, unlike the $M$
considered above, has $2^{g- 1}$ null rows corresponding to those
$\a$ for which $[\a,\b]$ is odd. To use an argument similar to the
one in the previous section relating $\P grTh$ and $\P Th_2$, we
need

\begin{lm}
For $g\geq 3$ the matrix $M^\b$ has rank at least three for all
$\tau$.
\end{lm}
\begin{proof}
Fix an irreducible period matrix $\tau$ and consider the abelian
variety $X:=X_{\tau}$. Let $\cal L$ denote the symmetric line
bundle inducing a principal polarization on $X$, for which $\tt 0
0 (\tau,z)$ is the basis for sections. We denote by $\Theta$ the
associated divisor. For any $x\in X$, let $t_x:X\to X$ be the
traslation by $x$. Setting $x:=\b/4$, we consider the line bundle
$\cal N:=t_{x}^{*}\cal L^2$. A basis for the sections of $\cal N$
are the theta functions $\tt \a\b (2\tau,2z)$. As a consequence of
(\ref{tT}),  we see that for any $x\in X$, all
$$
\tt 0\b(\tau,z+x)\tt 00(\tau,z-x)
$$
are linear combinations of $\tt\a\b(2\tau,2z)$. We recall that if
$\tau$ is irreducibile, the Gauss  maps $G_x:t_{x}^*\Theta\to
\P^{g-1}$ are dominant.

Suppose $M^\b$ is not of maximal rank. Then for some $\lambda$ and
$c$'s we have
$$
\sum_{i\leq j} c_{ij}\p_{\tau_{ij}} \tt \a\b (2\tau,0)=\lambda \tt
\a\b (2\tau,0).
$$
Thus we have
$$
\sum\limits_{i\leq j} c_{ij}\p_{z_{i}} \p_{z_{j}}\left.\left(\tt
0\b (\tau, z-x)\tt 00(\tau,z+x)\right)\right|_{z=0}
$$
$$
=\lambda\tt 0\b (\tau,x)\tt 00 (\tau,x)
$$
If $x$ does not belong to $\Theta\cup t_{\b/2}^*\Theta$, the
coefficient of $\lambda$ in the above relation is not zero. Vice
versa, assuming that $x\in \Theta\cap t_{\b/2}^*\Theta$, we have
$$
\sum\limits_{i\leq j} c_{ij}\p_{z_{i}}\tt 0\b (\tau,x)
\p_{z_{j}}\tt 00(\tau,x)+\p_{z_{j}} \tt 0\b (\tau,x)
\p_{z_{i}}\tt 00(\tau,x)=0.
$$
Now, since $\tau$ is irreducible, the singular locus of $\Theta$
has codimension at least two in $\Theta$, see \cite{el}, so there
exists such an $x$ in the smooth part of $\Theta$ and
$t_{\b/2}^*\Theta$. Moreover, for $g\ge 3$ by a linear
transformation we can find $x_1$ and $x_2$ such that
$$
\p_{z_{j}}\tt 00(\tau,x_1)=\de_{j}^1{\rm\ and\ }\p_{z_{j}}\tt
00(\tau,x_2) =\de_{j}^2,
$$
where $\de_i^j$ is Kr\"onecker's delta, and either $\p_{z_{j}}\tt
0\b(\tau,x_2)$  is not proportional to $\p_{z_{j}}\tt
00(\tau,x_1)$ or $\p_{z_{j}}\tt 0\b(\tau,x_1)$ is not proportional
to  $\p_{z_{j}}\tt 00(\tau, x_2)$. These properties impose three
linearly independent conditions on the coefficients
$c_{ij},\,\lambda$. Hence, for irreducible $\tau$ the matrix
$M^{\b}$ has at least three linearly independent columns. If the
point $\tau$ is reducible,  we can use  similar facts about the
Gauss map,  do directly the genus three  case and  use  some
inductive argument to finish the proof.
\end{proof}

Thus we get a generalization of proposition \ref{getTh2} in the
same way:
\begin{prop}
Given $\P grTh(\tau)$, the equations (\ref{syseq}) for $\b\neq 0$
have a unique projective solution for $\lbrace\tt\a\b(2\tau,0)
\rbrace_{\rm all\ \a}$.
\end{prop}

From the above proposition and formula (\ref{tT2}) it follows that
\begin{cor} All products of the type
$$
\tt\e{\b+\de}(\tau,0)\tt\e\de (\tau,0)
$$
are determined by $\P grTh$ uniquely up to a multiplicative
constant $t_{[0, \b]}$.
\end{cor}
Obviously the  above statement is true for every point $\tau$,
thus, using the $\Gamma_g$-equivariance of the map $\P grTh$ and
observing that the homogenuos action of $\Gamma_g$ on the  set of
characteristics is transitive on the set of characteristics
different from $[0,0]$, by using the above corollary stated for
points $\gamma\tau$ for all $\gamma$ we prove

\begin{cor}
\label{cor}
All products of the type
$$
\tt\a\b(\tau,0)\tt{\a+\e}{\b+\de}(\tau,0)
$$
are determined by $\P grTh(\tau)$ uniquely up to a multiplicative
constant $t_{\e, \de}$.
\end{cor}

Now we are able to prove our main theorem.

\begin{proof}[Proof of theorem 1]
We will show that generically  $\P grTh(\tau)$  determines $\P
Th(\tau)$ uniquely. Assume that there are two points $\tau$ and
$\tau'$ for which $\P grTh(\tau)=\P grTh(\tau').$ Since $\P
Th_2(\tau)=\P Th_2(\tau')$ by proposition \ref{getTh2},
$$
\tt\a\b(\tau,0)^2=c^2\tt\a\b(\tau',0)^2\qquad \forall [\a,\b],
$$
where $c$ is a constant independent of $\a,\b$. Hence
$$
\tt\a\b(\tau,0)=c s_{\a,\b}\tt\a\b(\tau',0),
$$
where $s_{\a, \b}$ is a sign depending on $\a,\b$. Replacing
$\tau$ by $\gamma\tau$ and $c$ by $-c$ if necessary, we can assume
$\tt 0 0 (\tau,0)\neq 0$ and $s_{0, 0}=1$. Now from the two
previous corollaries it follows that
\begin{equation}
\label{ok} t_{\a, \b}  = t_{0,0} s_{\a  ,\b },\qquad
s_{\a,\b}s_{\e,\de}=s_{\a+\e,\b+\de}\qquad \forall\a, \e,\de,\b
\end{equation}
whenever
\begin{equation}
\label{okif}
\tt\a\b(\tau,0)\cdot\tt\e\de(\tau,0)\cdot\tt{\a+\e}{\b+\de}(\tau,0)\ne0.
\end{equation}
\smallskip

We would like to prove that  all $s_{\a  ,\b }=1$. Using
(\ref{tT}) we can easily check that for all $[\e, \de]$, among the
products $\tt\a\b(\tau,0)\tt{\a+\e}{\b+\de}(\tau,0)$ appearing in
Corollary \ref{cor} there is at least one different from $0$.
Hence there is a linear basis
$[\e_1,\de_1],\dots,[\e_{2g},\de_{2g}]$ for the set of
characteristics such that the associated theta constants do not
vanish at $\tau$. We denote by $A_1(\tau)$ this
set of characteristics. Now recall that an element $\gamma=\pmatrix a&b\\
c&d\endpmatrix\in \Gamma_g(4)$ acts on theta constants by the
character
\begin{equation}
\label{char}
\phi(\e,\,\de,\gamma)=(-1)^{{\rm diag(b)\cdot\e}\over
4}(-1)^{{\rm diag(c)\cdot\de}\over 4}
\end{equation}
by the theta transformation law (\ref{translaw}). We observe that
on one side this action is compatible with the statement of
Corollary \ref{cor}. Thus we can choose $\gamma\in\Gamma_g(4)$
such that
$$
\tt{\e_j}{\de_j}(\tau,0)=c\,\tt{\e_j}{\de_j}(\gamma\tau ',0),
$$
i.e. $s_{\e_j,\de_j }=1$ for all $j=1,\dots, 2g$.

\smallskip
Now let $A(\tau)$ be the set of characteristics whose associated
theta constants do not vanish at $\tau$. We denote by $A_2(\tau)$
the subset of $A(\tau)$ consisting of the elements of $A_1(\tau)$
and the sums of two elements of $A_1(\tau)$ that are still in
$A(\tau)$. Since the condition (\ref{okif}) holds for
characteristics in $A_2(\tau)$, (\ref{ok}) is satisfied and thus
$s_{\e,\de }=1$ since this is true for $A_1(\tau)$.

Furthermore, let us denote by $A_3(\tau)$ the subset of $A$ whose
elements are either in $A_2(\tau)$ or sums of two elements of
$A_2(\tau)$. Iterating the process, we get a certain subset
$B(\tau)$. We would like to have $B(\tau)=A(\tau)$. Obviously this
holds and we have no trouble if we assume the non-vanishing of all
theta constants with characteristic $[\e,0]$ and $[0,\de ]$: then
whenever $\tt\e\de(\tau,0)\ne0$, we can use
$$
\tt\e0(\tau,0)\tt0\de(\tau,0)\tt\e\de(\tau,0)\ne0
$$
for the condition (\ref{okif}) and by formula (\ref{ok}) we are
done. So if we restrict to the open set determined by this
non-vanishing, we have $\P Th(\tau)=\P Th(\gamma \tau')$ for some
$\gamma\in\Gamma_g$, and thus from the injectivity of $\P Th$
follows that $\tau$ and $\gamma\tau'$ are
$\Gamma_g(4,8)$--conjugate, so that also $\P grTh(\tau)=\P
grTh(\gamma \tau')$. Now it is left to show that
$\gamma\in\Gamma_g(4,8)$.

\smallskip
Assume the contrary: $\gamma\in\Gamma_g(4)\setminus\Gamma_g(4,8)$.
We claim that then $\gamma$ acts non-trivially on $\P grTh(\tau)$,
so that we would have $\P grTh(\tau)=\P grTh(\gamma\tau')\not=\P
grTh(\tau')$, which is a contradiction. Indeed, $\gamma$ acts on
each gradient by multiplication by a sign $\phi(\e,\de,\gamma)$.
Consider all odd $[\e,\de]$ such that $\phi(\e,\de, \gamma)=-1$:
if for at least one of those and one of the remaining (otherwise
there is multiplication by $-1$) the associated gradient $\vec{\rm
grad}_z \tt\e\de$ is not the zero vector, we are done. Up to
conjugating by some element of $\Gamma_g$, we can assume that the
$b$ in $\gamma$ is such that ${\rm diag}(b)\equiv 0{\rm\ mod\ }8$.
Since all level subgroups are normal, a conjugate of $\gamma$ lies
exactly in the level subgroups in which $\gamma$ lies.

From formula (\ref{char}) we see that if ${\rm diag}\,
c\cdot\de\equiv 4{\rm\ mod\ }8$ for some fixed $\de$ (resp.
$\equiv 0$), then $\phi(\e, \de,\gamma)=-1$ (resp. $=1$) for all
$\e$. But since the map $X_\tau\to \P^{2^g-1}$ defined by
$z\to\tt\e0(\tau,z)$ is of maximal rank at the point $\de/2$
(recall that $\tt\e0(\tau,z+\de/2)=\tt\e\de(\tau,z)$), all
gradients at zero of odd theta functions of the form $\tt\e\de$
cannot vanish simultaneously for all $\e$.
\end{proof}

\begin{rem}
As a consequence of  the above  proof, if we have
$B(\tau)=A(\tau)$ for all $\tau$, the map $\P grTh$ is injective.

To the best of our knowledge this is always true. In particular it
easy to check that even when $\tau$ is $\Gamma_g$-conjugate to a
diagonal matrix, i.e. is the period matrix of a product of
elliptic curves and thus has the maximal possible number of
vanishing theta constants.

Thus we can say that $\P grTh$ is injective, if for  all $\tau$
the  corresponding subset $A(\tau)$ contains $A(\tilde\tau)$, with
$\tilde\tau$ being $\Gamma_g$-conjugate to a diagonal matrix. For
all examples of abelian varieties that we know the set of
characteristics for which the associated theta constants vanish is
always contained in such a set for a diagonal period matrix, up to
conjugation. Thus in all the examples that we know we do have
$A(\tau)\supset A(\tilde\tau)$ and the map $\P grTh$ is injective
at $\tau$.
\end{rem}

Since for $g=3,\, 4$ the combinatorics of the possible vanishing
of theta constants is known well, and the worst cases are the
reducible and the hyperelliptic, which can be treated by hand, we
have
\begin{cor}
For $g=3,\, 4$ the map $\P grTh$ is injective on $\A_g^{4,8}$.
\end{cor}

We also prove the injectivity of $\P grTh$ on generic Jacobians:
notice that it does not directly follow from theorem 1 a priori as
the Jacobians may not be the ``generic'' abelian varieties.

\begin{proof}[Proof of corollary 2]
From the proof of theorem 1 we see that the map $\P grTh$ is
injective at $\tau$ if none of the theta constants at $\tau$
vanish. However, it is known classically from the works of Riemann
that no theta constant vanishes identically on $\J_g$
Since $\J_g$ is irreducible, the subset of $\J_g$ where no theta
constants vanish is Zariski open, and there the injectivity of $\P
grTh$ holds.
\end{proof}

\section{Injectivity on the tangent space}
Now let us show that $\P grTh$ is smooth. Let us study the
situation in general terms first.
\begin{lm}
Suppose $\P f:X\to Gr(k,N)$ is an analytic map of a complex
variety to the grassmannian, locally near some $p\in X$ given by
$f:x_1\ldots x_M\to f_1(x), \ldots f_k(x)$, where $x_i$ are the
local coordinates near $p$, and each $f_i$ is a vector in $\C^N$.
Then $\dd \P f|_p$ is injective if and only if for all
$v\in\C^M\setminus\lbrace 0\rbrace$ at least for one $I$ the
vector $\p_v f_I(0)$ is linearly independent with $(f_1(0),\ldots,
f_k(0))$. In particular it is injective if the vectors
$(f_1(0),\ldots, f_k(0),\p_{x_1}f_I(0),\ldots, \p_{x_M}f_I(0))$
are linearly independent for some $I$.
\label{diffinj}
\end{lm}
\begin{proof}
Indeed let us consider the linearization of $f$ near $p$; for
$v\in \C^M$ infinitesimally small we have
$f(v)=(f_1(0)+\p_vf_1(0),\ldots, f_N(0)+\p_vf_N(0))$. $\dd \P
f|_p$ is non-degenerate iff it does not map any tangent vector $v$
to zero, i.e. if $f(v)$ represents a point in the grassmannian
different from $f(0)$ for all $v$. This is equivalent to saying
that for any $v$ at least one of the vectors making up $f(v)$ does
not lie in $\P f(0)$, i.e. that at least one $\p_vf_I(0)$ is
linearly independent with $(f_1(0),\ldots, f_k(0))$.

Now if $(f_1(0),\ldots, f_k(0),\p_{x_1}f_I(0),\ldots, \p_{x_M}
f_I(0))$ are linearly independent for some $I$, then it implies
that any linear combination
$$
\p_vf_I(0)=\sum v_j\p_{x_j}f_i(0)
$$
is linearly independent with $(f_1(0),\ldots, f_k(0))$ for all
$v$.
\end{proof}

\begin{proof}[Proof of theorem 3: injectivity of $\dd \P grTh$ on
tangent spaces] We use the above lemma for $f=grTh$, $k=g$,
$M=g(g+1)/2$ and $N=2^{g-1}(2^g-1)$. In \cite{sm2} it is stated
(Theorem 2b) that $\P grTh$ is an immersion away from
$\Gamma_g$-conjugates of points that are reducible as a product of
a one-dimensional and a $g-1$-dimensional abelian variety. The
proof there proceeds by showing that the rank of the matrix
$(grTh(\tau),\p_{\tau_{ij}}grTh(\tau))$, (i.e. of
$(f_i(0),\p_{x_j}f_i(0))_{\rm all\ i,j}$ in the notations of the
lemma \ref{diffinj}) is maximal exactly for points that are not
one-reducible. The maximality of this rank implies the maximality
of the rank of any submatrix, in particular the one for which we
need linear independence to apply lemma \ref{diffinj}. However,
the converse, implicitly assumed in \cite{sm2}, is in fact not
true: it would be requiring in addition to the lemma that
$(\p_{x_j}f_i(0))$ are linearly independent, which is not
necessary.

Thus the argument in \cite{sm2} shows that $\P grTh$ is an
immersion away from $\H_1\times \H_{g-1}$, and we only need to
deal with $\tau=\pmatrix \lambda&0\\ 0&\tau'\endpmatrix$. In this
case write $[\e,\de]=[\e_1\,\e',\de_1\,\de']$, and let indices
$i,j,I$ always be greater than one. Then
$$
\p_{z_1}\tt\e\de(\tau)=\p_{z_1}\tt{\e_1}{\de_1}(\lambda)\tt{\e'}{\de'}
(\tau');
\ \p_{z_I}\tt\e\de(\tau)=\tt{\e_1}{\de_1}(\lambda)\p_{z_I}\tt{\e'}{\de'}
(\tau')
$$
Using the heat equation, for $\dd\, \p_{z_1}\tt\e\de|_\tau$ we have
$$
\p_{\tau_{ij}}\p_{z_1}\tt\e\de=\p_{z_1}\tt{\e_1}{\de_1}\p_{\tau_{ij}}
\tt{\e'}{\de'};\ \p_{\tau_{11}}\p_{z_1}\tt\e\de=\p_{z_1}^3
\tt{\e_1}{\de_1}\tt{\e'}{\de'}
$$
$$
\p_{\tau_{1i}}\p_{z_1}\tt\e\de=\p_{\tau_{11}}\tt{\e_1}{\de_1}\p_{z_i}
\tt{\e'}{\de'}.
$$
What are the possible linear relations among these vectors?
Arrange the vectors into a matrix and split this matrix into two
corresponding to whether $[\e',\de']$ is odd or even. We notice
that all derivatives above are non-zero only for one parity of
$[\e',\de']$ (i.e. every column of the matrix has non-zero
elements only in one of the two submatrices), and thus the matrix
of vectors is in $2\times 2$ submatrix block form, and we can
compute the rank by adding the ranks of the blocks.

For $[\e',\de']$ odd, the only non-zero elements of the
corresponding row are $\p_{z_i}\tt\e\de$ and
$\p_{\tau_{1i}}\p_{z_1}\tt\e\de$, which are independent because
all $\p_{z_i}\tt{\e'} {\de'}$ are non-collinear (as $grTh(\tau')$)
and the matrix $(\tt{\e_1}{\de_1},\p_{\tau_{11}}
\tt{\e_1}{\de_1})$ has maximal rank, two.

On the other hand, if $[\e',\de']$ is even, then
$\p_{z_1}\tt\e\de$ and $\p_{\tau_{ij}}\p_{z_1}\tt\e\de$ are
independent because the matrix $(\tt{\e'}{\de'},
\p_{\tau_{ij}}\tt{\e'}{\de'})$ has maximal rank, while
$\p_{z_1}^3\tt\e\de$ is proportional to $\p_{z_1}\tt\e\de$.

Thus using the above lemma we see that if $v\not\in\C\p_{\tau_{
11}}$, the vectors $(f_1(0),f_i(0), \p_vf_1(0))$ are linearly
independent. By the lemma, to prove the injectivity of $\dd \P f$
we then need to show that for $v=\p_{\tau_{11}}$ the vectors
$(f_1(0),f_i(0),\p_vf_I(0))$ are linearly independent. Indeed we
compute
$$
\p_{\tau_{11}}\p_{z_I}\tt\e\de(\tau)=
\p_{\tau_{11}}\tt{\e_1}{\de_1}(\lambda) \p_{z_I}\tt{\e'}{\de'},
$$
which is linearly independent with $\p_{z_I}\tt\e\de$ because
$(\tt{\e_1}{\de_1},\p_{\tau_{11}}\tt{\e_1}{\de_1})$ has maximal
rank, and with $\p_{z_1}\tt\e\de$ because one is zero for
$[\e',\de']$ odd, and the other --- for $[\e',\de']$ even.
\end{proof}

\section*{Acknowledgements}
We are very grateful to Lucia Caporaso and Edoardo Sernesi for
bringing the subject to our attention and explaining to us their
work on theta-hyperplanes, which inspired this article. We would
also like to thank the organizers of the Complex Analysis meeting
in Oberwolfach in August 2002 for bringing professor Caporaso and
the first author to the same place at the same time, and thus
making this work begin.


\begin{thebibliography}{Fr885}
\bibitem[Ar872]{aronh}{{\it S. Aronhold,} Sur les vingt-huit tangentes
doubles d'une courbe du quatri\`eme degr\'e. Nouv. Ann. (2)  XI.
438--443. (1872).}
\bibitem[Ca01]{cap1}{{\it L. Caporaso,} On modular properties of odd
theta-characteristics, Advances in algebraic geometry motivated by
physics, 101--114, Contemp. Math., 276, Amer. Math. Soc.,
Providence, RI, 2001}
\bibitem[CS00]{cap2}{{\it L. Caporaso, E. Sernesi,} Recovering plane
curves from their bitangents, J. Algebraic Geom. {\bf 12} (2003),
no. 2, 225--244.}
\bibitem[CS02]{cap4}{{\it L. Caporaso, E. Sernesi,} Characterizing curves by
their odd theta-characteristics, math.AG/0204164}
\bibitem[EL97]{el}{{\it L. Ein, R. Lazarsfeld,} Singularities of theta divisors and
the birational geometry of irregular varieties, J. Amer. Math.
Soc. {\bf 10} (1997), no. 1, 243--258.}
\bibitem[Fa79]{fay}{{\it J. Fay,} On the Riemann-Jacobi formula,
Nachr. Akad. Wiss. Gottingen Math.-Phys. Kl. II 1979, no. 5, 61--73}
\bibitem[Fr885]{fr885}{{\it F. Frobenius,} Uber die constanten Factoren der
Thetarheinen, J. Reine Angew. Math. {\bf 98} (1885) 244--265}
\bibitem[Ig72]{ig1}{{\it J.-I. Igusa,} Theta functions. Die Grundlehren
der mathematischen Wissenschaften, Band 194. Springer-Verlag, New
York-Heidelberg, 1972.}
\bibitem[Ig80]{ig2}{{\it J.-I. Igusa,} On Jacobi's derivative formula
and its generalizations, Amer. J. Math. {\bf 102} (1980), no. 2,
409--446}
\bibitem[SM83]{sm3}{{\it R. Salvati Manni,} On the nonidentically zero
Nullwerte of Jacobians of theta functions with odd
characteristics, Adv. in Math. {\bf 47} (1983), no. 1, 88--104.}
\bibitem[SM94]{sm1}{{\it R. Salvati Manni,}  Modular varieties with level
$2$ theta structure, Amer. J. Math. {\bf 116} (1994), no. 6,
1489--1511.}
\bibitem[SM96]{sm2}{{\it R. Salvati Manni,} On the differential of
applications defined on the moduli space of p.p.a.v. with level
theta structure, Math. Z. {\bf 221}, 231-241 (1996)}
\bibitem[Sa83]{sasaki}{{\it R. Sasaki,} Modular forms vanishing at the
reducible points of the Siegel upper-half space, J. Reine Angew.
Math. {\bf 345} (1983), 111--121}
\end{thebibliography}
\end{document}